\documentclass[12pt, final]{amsproc}
\usepackage[cp850]{inputenc}
\usepackage{amssymb}
\usepackage[mathscr]{eucal}
\usepackage{amscd}
\usepackage{amsmath}
\usepackage{setspace}
\usepackage{tikz}
\def\Ker{\operatorname{Ker}}

\def\dim{\operatorname{dim}}
\def\supp{\operatorname{supp}}

\theoremstyle{plain}
\newtheorem{theorem}{Theorem}
\newtheorem{proposition}{Proposition}
\newtheorem{lemma}{Lemma}
\newtheorem{corollary}{Corollary}

\theoremstyle{remark}

\newtheorem{example}{Example}

\title{A UNIVERSAL FORMULA FOR DERIVATIONS OPERATORS AND APPLICATIONS}

\address{Departamento de Matem\'aticas, Universidad de Extremadura, Avenida de Elvas s/n, 06011 Badajoz, Spain.}
\email{jesus@unex.es}
\author{Jes\'us Su\'arez de la Fuente}


\thanks{The author was supported by no institution.}
\begin{document}
\begin{abstract} We give a universal formula describing derivation operators on a Hilbert space for a large class of interpolation methods. It is based on a simple new technique on ``critical points" where all the derivations attain the maximum. We deduce from this a version of Kalton uniqueness theorem for such methods, in particular, for the real method. As an application of our ideas is the construction of a weak Hilbert space induced by the real $J$-method. Previously, such space was only known arising from the complex method. To complete the picture, we show, using a breakthrough of Johnson and Szankowski, nontrivial derivations whose values on the critical points grow to infinity as slowly as we wish.
\end{abstract}
\maketitle

\section{Introduction}
In the 90's a very special class of non-linear maps, the so called derivation operators, which are a very special class of nonlinear maps, were intensively studied in the setting of the complex interpolation method. Nigel Kalton delivered the \textit{coup de gr\^ace} after the publication of his monumental \textit{``Differentials of complex interpolation processes for K\"othe function spaces"} cf. \cite{ka}. The idea of derivation will lead the reader, sooner or later, to a nowadays standard way to construct the Kalton-Peck space \cite{KaPe}. This is done with the so-called Kalton-Peck derivation $$\Omega(a)=a\cdot \log a.$$

One of the ideas surrounding this paper is that derivation operators do not belong exclusively to the complex method. An easy way to notice this is the general description of interpolation methods given by Cwikel, Kalton, Milman and Rochberg \cite{CKMR}. While they outline a portrait of many interpolation methods under the same schema (including the real interpolation method), this schema follows clearly a ``complex roadmap". Thus, in order to construct derivation operators we only need to follow the signals to construct derivations, for example, for the real method. Indeed, the approach of \cite{CKMR} sinks its roots also in the discrete version of the Lions and Peetre  \textit{espaces de moyennes} cf. \cite{LiPe} (see also Subsection \ref{kaltonpeque}). This is, as it is well known, a real method.

The derivation operators can be used to construct unexpected examples of spaces as the Kalton-Peck space which force us to broaden our view. For example, using derivations, the author constructed a new example of a weak Hilbert space which is an important class of spaces introduced by Pisier in \cite{Pi}. This example, denoted $Z(T^2)$ (cf. \cite{JS}), appears as the derivation of the complex interpolation formula $(T^2, (T^2)^*)_{\frac{1}{2}}=\ell_2$ (see \cite{CoS}), where $T^2$ denotes the $2$-convexification of the Tsirelson space. Following Cwikel, Milman and Rochberg \cite{CMR}, the resulting derivation operator, say $\Omega$, represents the tangent vector at the instant $1/2$ to a ``Calder\'on curve" joining $T^2$ and its dual $(T^2)^*$, see also Semmes \cite{Se}. To avoid trivialities when the derivation operator is zero, Kalton provided us with the so-called Kalton uniqueness theorem \cite[Theorem 7.6.]{ka}: 
\begin{theorem}
Let $X=(X_0,X_1)$ and $Y=(Y_0,Y_1)$ be compatible couples of K\"othe sequence spaces such that $(X_0,X_1)_{\frac{1}{2}}=\ell_2=(Y_0,Y_1)_{\frac{1}{2}}$ and let $\Omega_X$ and $\Omega_Y$ be the corresponding derivation operators. If $\Omega_X$ and $\Omega_Y$ are equivalent, then $X_0=Y_0$ and $X_1=Y_1$ with an equivalent norm.
\end{theorem}
To explain clearly the meaning of this result pick $Y_0=Y_1=\ell_2$ so that $\Omega_Y=0$ (since the curve is constant the derivative must be zero). Thus, if $\Omega_X$ is equivalent to $\Omega_Y=0$ then $X_0=X_1=\ell_2$. Let us reduce the claim to its foundations: if the derivative is zero, then the variable is constant. The reader will agree that, at least stated in this way, the claim seems not to be as something exclusive of the complex method. However, the sad point to deal with derivation operators is that there is not an analogue of Kalton's result for all interpolation methods. There is a result for the real method (\cite[Theorem 5.17]{CJMR}) but the proof is long and works only for that method.

 We give a simple and unified proof for most of the interpolation methods fitting in the schema of Cwikel, Kalton, Milman and Rochberg \cite{CKMR}. The first step is the computation of the derivations for all these interpolation methods. Thus, let us fix a couple $(B,B^*)$ of sequence spaces and an interpolation method $\mathcal F$ so that $\mathcal F (B,B^*)=\ell_2$. While we cannot compute the corresponding derivation operator $\Omega$ in all the points, we shall restrict ourselves to $\mathbb C^n$ and calculate it at the points $a_n\neq 0$ where the maximum of the derivation operator is attained by compactness. These ``critical points" $a_n$ are independent of the chosen method $\mathcal F$ as they only obey to the equivalence constant, say $\kappa(n)$, of $\ell_2$ with $B$ on $\mathbb C^n$. To put it tersely, once the couple is fixed, all the derivations arising from our large class of interpolation methods are $$\Omega(a_n)=a_n \cdot \log \kappa(n).$$ 
The factor $\log \kappa(n)$ has a geometric flavor since it can be read as the ``Calder\'on distance" introduced by Semmes (cf. \cite{Se}, see also the discussion in Subsection \ref{geometric}) between the norms of $B$ and the Hilbert on $\mathbb C^n$. 

As a matter of fact, the formula shows that the derivations, on the critical points, behave asymptotically in the same way. The rate of growth is independent of the interpolation method and it only depends on the Calder\'on distance which is determined just by the fixed couple $(B,B^*)$ and not by the interpolation method itself. We would like to recall that the form of $\Omega$ is, in general, unknown for a couple $(B,B^*)$ except for a very few cases as the Kalton-Peck derivation. However, the value $\kappa(n)$ is easy to estimate for many spaces.

Kalton uniqueness theorem follows now readily for $\mathcal F (B,B^*)=\ell_2$. If the derivation is bounded then so does $\log \kappa(n)$ for every $n\in \mathbb N$. This implies that $\kappa(n)$ and $\kappa(n)^{-1}$ are bounded so that $B=B^*=\ell_2$.
We can easily deduce now that the derivation corresponding to the interpolation space $(T^2, (T^2)^*)_{\frac{1}{2},2;J}=\ell_2$ by the real $J$-method is not trivial. 

We introduce the method of critical points in Section \ref{sec4} with two enlightening examples: the real and the complex method. The universal formula and our version of Kalton uniqueness theorem is given in full detail in Section \ref{sec5}.  To convince us that the ``real" derived space of $(T^2, (T^2)^*)_{\frac{1}{2},2;J}$ is still a weak Hilbert space (as in the complex case), we have to develope some probabilistic estimates for the real method as the ones given in \cite{JS} for the complex method. We believe that these estimates are of independent interest and can be found at Subsection \ref{sec3}. After this, we may give a proof that the derived space $d(T^2, (T^2)^*)_{\frac{1}{2},2;J}$ is also a weak Hilbert space.

 Until here we have computed derivations for a fixed couple. To conclude, in the last subsection, we go in the opposite direction. We show, using important results of Johnson and Szankowski \cite{JoSa}, how given a sequence $1\leq \delta_n\to \infty$, we may produce derivation operators $\Omega$ whose ``maximum values", say $\|\Omega(a_n)\|$, on the corresponding critical points $(a_n)$, are bounded by $\delta_n$, up to universal constant. So, in particular, we find a rich variety of such maps. This closes our circle of ideas.


\section{Preliminaries}
\subsection{Notation}
We use standard notation for Banach spaces as provided in the book of Lindenstraus and Tzafriri \cite{LT} or in Albiac and Kalton \cite{AK}. The word space is reserved for Banach space. We work with sequence spaces so if no confusion arises $(e_j)_{j=1}^{\infty}$ will denote the natural basis. Given a sequence $(x_j)_{j=1}^{\infty}$ of a space, we denote $[x_j]_{j=1}^{\infty}$ its closed linear span. We shall use the compact notation $[1,n]:=\{1,...,n\}$ and $\lfloor t \rfloor$ will denote the floor function, where $t\in \mathbb R$. Also, given two sequences of real numbers $(a_j)_{j=1}^{\infty}$ and $(b_j)_{j=1}^{\infty}$, we write $a_j \sim b_j$ (respectively $a_j \lesssim b_j$) if there is $c>0$ such that $c^{-1}\cdot b_j\leq a_j\leq c \cdot b_j$ (respectively $a_j\leq c\cdot b_j)$ for all $j\in \mathbb N$. Recall that for a space $X$ and a fixed $1<p\leq 2$, the number $a_{m,p}(X)$ is defined to be the least constant $a$ such that
$$\mathbb E\left \|\sum_{j=1}^m \varepsilon_j x_j\right\|\leq a \left ( \sum_{j=1}^m \| x_j\|^p\right)^{1/p},$$
for all $x_1,...,x_m\in X$ and where the average is taken over all $(\varepsilon_j)_{j=1}^m\in \{-1,1\}^m$. Thus, the space $X$ has \textit{type $p$} if $a_p(X):=\sup_{m\in \mathbb N}a_{m,p}(X)<\infty$.
To finish, we say $X$ is a weak Hilbert space if there is $0<\delta_0<1$ and a constant $C$ with the following property: every finite dimensional subspace $E$ of $X$ contains a subspace $F\subseteq E$ with $\dim F\geq \delta_0 \dim E$ such that $d_F\leq C$ and there is a projection $P:X\to F$ with $\|P\|\leq C$.
The definition above is not the original one but is chosen out among the many equivalent characterizations given by Pisier \cite[Theorem 12.2.(iii)]{Pib}. Recall that $T^2$, the $2$-convexification of the Tsirelson space, is the most important example of weak Hilbert space and the reader may found in \cite{CS} or \cite{Pib} a comprehensive study. 
\subsection{Interpolation affairs} We follow closely the approach of Cwikel, Kalton, Milman and Rochberg \cite{CKMR}.
\subsubsection{The pseudolattice and the Calder\'on space}
Let $\mathbf{Ban}$ be the class of all Banach spaces over the complex numbers. A mapping $\mathcal X:\mathbf{Ban} \to \mathbf{Ban}$ will be called a pseudolattice if
\begin{enumerate}
\item[(i)] for each $B\in \mathbf{Ban}$ the space $\mathcal X(B)$ consists of $B$-valued sequences $\{b_n\}_{n\in \mathbb Z}$ and if
\item[(ii)] whenever $A$ is a closed subspace of $B$ it follows that $\mathcal X(A)$ is a closed subspace of $\mathcal X(B)$ and if
\item[(iii)] there exists a positive constant $C>0$ such that, for all $A,B\in \mathbf{Ban}$, and all bounded linear operators $T:A\to B$ and every sequence $\{a_n\}_{n\in \mathbb Z}\in \mathcal X(A)$, the sequence $\{T(a_n)\}_{n\in \mathbb Z}\in \mathcal X(B)$ satisfies the estimate $$\| \{T(a_n)\}_{n\in \mathbb Z} \|_{\mathcal X(B)}\leq C(\mathcal X)\|T\|_{A\to B}\|\{a_n\}_{n\in \mathbb Z}\|_{\mathcal X(A)} .$$
\end{enumerate}
As an example, let $X$ denote a Banach lattice of real valued functions defined over $\mathbb Z$; the reader may keep in mind $X=\ell_2$.  We will use the notation $\mathcal X=X$ to mean that, for each $B\in  \mathbf{Ban}$, $\mathcal X(B)$ is the space, also denoted $X(B)$, of all $B$-valued sequences $\{b_n\}_{n\in \mathbb Z}$ such that $\{\|b_n\|_X\}_{n\in \mathbb Z}\in X$ endowed with the norm $$\|\{b_n\}_{n\in \mathbb Z}\|_{X(B)}=\|\{\|b_n\|_B\}_{n\in \mathbb Z}\|_X.$$
Fix now a pair of pseudolattices  $\mathbf {X}=\{ \mathcal X_0, \mathcal X_1\}$. Given a compatible couple of Banach spaces $B=(B_0,B_1)$ we define $$\mathcal J(\mathbf {X},B)$$ to be the space of all $(B_0\cap B_1)$-valued sequences $\{b_n\}_{n\in \mathbb Z}$ for which the sequence $\{e^{jn}b_n\}_{n\in \mathbb Z}\in \mathcal X_j(B_j) $ for $j=0,1$. This space is normed by 
$$\| \{b_n\}_{n\in \mathbb Z}\|_{\mathcal {J} (\mathbf {X},B)}=\max_{j=0,1} \| \{e^{jn}b_n\}\|_{\mathcal X_j(B_j)}.$$
\subsubsection{Nontrivial and Laurent compatible pseudolattices} The pseudolattice pair $\mathbf {X}$ is \textit{nontrivial} if, for the special one-dimensional Banach pair $(\mathbb C, \mathbb C)$ and each $z\in \mathbb A=\{z\in \mathbb C:1<|z|<e \}$ (the open annulus) there exists $\{b_n\}_{n\in \mathbb Z}\in \mathcal {J} (\mathbf {X},(\mathbb C, \mathbb C)) $ such that the series $\sum_{n\in \mathbb Z}z^nb_n$ converges to a nonzero number. The pseudolattice pair $\mathbf {X}$ is \textit{Laurent compatible} if it is nontrivial and for every $z\in \mathbb A$ the Laurent series $$\sum_{n\in \mathbb Z} z^nb_n$$ converges absolutely with respect to the norm of $B_0+B_1$. Therefore the sum of this series is an analytic function of $z$ in $\mathbb A$ and can be differentiated term-by-term. The series for its derivative $$\sum_{n\in \mathbb Z} nz^{n-1}b_n$$ also converge absolutely in $B_0+B_1$. See both claims at \cite[page 248]{CKMR}.
\subsubsection{The interpolation space and the evaluation map}
Let $B=(B_0,B_1)$ be a compatible pair, then for each $0<\theta<1$, we define the \textit{interpolation space} $B_{\mathbf {X}, \theta}$ to consist of all elements of the form $$b=\sum_{n\in \mathbb Z}e^{\theta n}b_n$$ with $\{b_n\}_{n\in \mathbb Z}\in \mathcal J(\mathbf {X},B)$, endowed with the natural quotient norm:
$$\| b \|_{B_{\mathbf {X}, \theta}}=\inf \left\lbrace \| \{b_n\}_{n\in \mathbb Z}\|_{\mathcal J(\mathbf {X},B)}: b= \sum_{n\in \mathbb Z}e^{\theta n}b_n   \right\rbrace .$$
According to the previous claims one may ``think" every element $\{b_n\}_{n\in \mathbb Z}\in \mathcal J(\mathbf {X},B)$ as the analytic map $\sum_{n\in \mathbb Z} z^nb_n$, where $z\in \mathbb A$, with all the precautions. Hence we shall informally write $$\{b_n\}_{n\in \mathbb Z}=\sum_{n\in \mathbb Z} z^nb_n.$$  Therefore, we have the natural \textit{evaluation map} $$\delta_{\theta}:\mathcal J(\mathbf {X},B) \longrightarrow B_{\mathbf {X}, \theta}$$ given by the rule 
$$\delta_{\theta}(\{b_n\}_{n\in \mathbb Z})=\sum_{n\in \mathbb Z} e^{\theta n} b_n.$$
Two examples to keep in mind are the complex and real method, respectively. The complex method corresponds to $\mathcal X_0=\mathcal X_1=FC(\mathbb T,\cdot)$ which denotes the Fourier transform of continuous functions over the torus $\mathbb T$ endowed with the sup-norm. The real method correspond to the choice $\mathcal X_j=\ell_{p_j}$ with $j=0,1$. Both claims are easy to find in \cite{CKMR}.
\subsubsection{The structure of module}
We say a pseudolattice pair $\mathbf X$ is an $\ell_{\infty}$-module if whenever $(B_0,B_1)$ are spaces with the same $1$-unconditional basis there is $C>0$ such that for every $\{b_n\}_{n\in \mathbb Z}\in \mathcal J(\mathbf {X},B)$ and $a\in \ell_{\infty}$, the following holds:
\begin{enumerate}
\item $\{a\cdot b_n\}_{n\in \mathbb Z}\in \mathcal J(\mathbf {X},B)$
\item $\|\{a\cdot b_n\}_{n\in \mathbb Z}\|_{\mathcal J(\mathbf {X},B)}\leq C\|a\|_{\infty}\|\{b_n\}_{n\in \mathbb Z}\|_{\mathcal J(\mathbf {X},B)}$
\end{enumerate}
An example of this is $\mathbf X=\{\ell_{p_0}, \ell_{p_1}\}$ because one trivially has $$\|a\cdot b_n\|_{B_j}\leq \|a\|_{\infty} \|b_n\|_{B_j},\;\;j=0,1.$$
And thus for $j=0,1$, $$\|\{a\cdot e^{jn} b_n\}_{n\in \mathbb Z}\|_{\ell_{p_j}(B_j)} \leq \|a\|_{\infty}\| \{e^{jn} b_n\}_{n\in \mathbb Z}\|_{\ell_{p_j}(B_j)}.$$

\subsubsection{The derivation operator and centralizers}
Let us fix a constant $C>1$. For a given $b\in B_{\mathbf {X}, \theta}$, pick $\{b_n\}_{n\in \mathbb Z}$ so that $\delta_{\theta}(\{b_n\}_{n\in \mathbb Z})=b$ and $$\|\{b_n\}_{n\in \mathbb Z}\|\leq C\|b\|.$$ Sometimes we refer to $\{b_n\}_{n\in \mathbb Z}$ as an \textit{extremal}. We write $$S(b)=\{b_n\}_{n\in \mathbb Z}$$ and say that $S$ is $C$-\textit{bounded selector} for the map $\delta_{\theta}$. The \textit{derivation operator} is given by the rule
\begin{equation}\label{deriva}
\Omega(b):=\delta_{\theta}'S(b)=\sum_{n\in \mathbb Z} ne^{\theta(n-1)}b_n.
\end{equation}
If $W$ is the ambient space of the couple $B$, then $$\Omega:B_{\mathbf {X}, \theta}\longrightarrow W.$$ We deal with sequences spaces so our ambient space will be $\omega$, the vector space of all complex scalar sequences. We will work under the condition that $\mathbf {X}$ \textit{admits differentiation},  \cite[Lemma 3.11]{CKMR} which is a rather technical condition. However, the reader may keep in mind that $\mathbf {X}$ admits differentiation if the shift operator is an isometry on $\mathcal X_j(B_j)$ for $j=0,1$, see \cite[Lemma 3.6.]{CKMR}. This will be our case: the shift operator is clearly an isometry on $\ell_{p_j}(B_j)$ for $j=0,1$, see \cite[Lemma 3.6.]{CKMR}.
%
Let $\mathbf X$ be a pseudolattice pair that is an $\ell_{\infty}$-module. For spaces $(B_0,B_1)$ with the same $1$-unconditional basis, the corresponding derivation operator $\Omega$ for $B_{\mathbf X, \theta}$ is a \textit{centralizer}, meaning that there is $C>0$ so that 
\begin{equation}\label{derivation}
\|\Omega(a\cdot b)-a\cdot \Omega(b)\|_{B_{\mathbf X, \theta}}\leq C\|a\|_{\infty}\|b\|_{B_{\mathbf X, \theta}},
\end{equation}
for all $a\in \ell_{\infty}$ and $b\in B_{\mathbf X, \theta}$. 
\subsubsection{Nontrivial and bounded derivations}
A derivation operator $\Omega:B_{\mathbf X, \theta}\to W$ is trivial if there is a linear (possibly unbounded) map $L:B_{\mathbf X, \theta}\to W$ such that $$\sup_{\|x\|\leq 1}\|\Omega(x)-L(x)\|<\infty.$$
Otherwise is \textit{nontrivial}. We also say two derivations are equivalent if the difference is a trivial derivation. For the particular case $L=0$, we say that the derivation $\Omega$ is \textit{bounded}. This is, if there is $C>0$ such that for every $b\in  B_{\mathbf X, \theta}$, we have that $\Omega(b)\in B_{\mathbf X, \theta}$ with $$\|\Omega(b)\|_{ B_{\mathbf X, \theta}}\leq C\|b\|_{ B_{\mathbf X, \theta}}.$$

\subsubsection{The derived space and the short exact sequence it generates}
Let $B$ be a compatible couple with ambient space $W$. The derived space $dB_{\mathbf X,\theta}$ is the set of couples $(x,y)\in W\times W$ for which the following quasi-norm  $$\|x-\Omega(y)\|_{B_{\mathbf X,\theta}}+\|y\|_{B_{\mathbf X,\theta}}$$
is finite, where $\Omega$ is the derivation operator (\ref{deriva}). If $B_{\mathbf X,\theta}$ contains no (uniform) copies of $\ell_1^n$ for every $n\in \mathbb N$ then such a quasi-norm is equivalent to a norm by a result of Kalton \cite{K4}. In our case, $B_{\mathbf X,\theta}=\ell_2$ so it is trivially satisfied. 
Recall that a \textit{short exact sequence} is a diagram of the form
\begin{equation}\label{ses}
\begin{CD} 0 @>>>Y@>j>>Z@>q>>X@>>>0
\end{CD}
\end{equation} where the morphisms are linear and continuous and such that the image of each arrow is the kernel of the next one. This condition implies that $Y$ is a subspace of $Z$ through $j$ and due to
the open mapping theorem we find that $X$ is isomorphic to $Z/j(Y)$. We usually refer to $Z$ as a \textit{twisted sum} of $Y$ and $X$. It is not hard to check that there is a short exact sequence 
\begin{equation*}
\begin{CD} 0 @>>>B_{\mathbf X,\theta}@>j>>d B_{\mathbf X,\theta}@>q>>B_{\mathbf X,\theta}@>>>0
\end{CD}
\end{equation*}
where $jx=(x,0)$ and $q(x,y)=y$. Let us observe that if $\Omega$ is bounded, then $$dB_{\mathbf X,\theta} \approx B_{\mathbf X,\theta}\oplus B_{\mathbf X,\theta},$$
where $\approx$ stands for isomorphic Banach spaces.
\subsubsection{Summary}
All the pseudolattice pairs $\mathbf X$ used throughout the paper will be assumed to be:
\begin{itemize}
\item Nontrivial.
\item Laurent compatible.
\item Admits differentation.
\end{itemize}
Our main examples are the pseudolattices corresponding to the real method that satisfies all the conditions above. The interpolated space by the $(\theta,q)$-real method corresponds to the choice $\mathcal X_0=\mathcal X_1=\ell_q$ or, equivalently, $\mathcal X_j=\ell_{q_j}$ for $j=0,1$ with $q^{-1}=(1-\theta)q_0^{-1}+\theta q_1^{-1}$, see the discussion in \cite[Paragraphs 2 and 4, page 251]{CKMR}. This is, both choices provide isomorphic Banach spaces. However, they give in general different derivation maps as we shall see in Subsection \ref{word}.
\subsection{Twisted Hilbert spaces}
 Let us assume that $Y\approx\ell_2\approx X$ in (\ref{ses}), then we simply say that $Z$ is a \textit{twisted Hilbert space}. The most important example of twisted Hilbert space is the so-called Kalton-Peck space $Z_2$ \cite{KaPe}. It appears from the complex method as the derived space of $(c_0,\ell_1)_{1/2}=\ell_2$. The corresponding derivation is the Kalton-Peck map, given as $$\Omega(b)=2\cdot\sum_{j=1}^{\infty} b_j\log \frac{|b_j|}{\|b\|}e_j.$$
We are also interested in the space $Z(T^2)$ which is the first example of a nontrivial twisted Hilbert space that is a weak Hilbert space \cite{JS}. It is defined as the derived space of $(T^2, (T^2)^*)_{1/2}=\ell_2$ by the complex method.

We finish recalling a couple of results that will be useful for us. For the first, let $(B,B^*)$ do share the same unconditional basis and let $\mathbf X$ be a pseudolattice pair that is an $\ell_{\infty}$-module. Let us assume that $(B,B^*)_{\mathbf X, 1/2}=\ell_2$, so that the corresponding derivation operator is a centralizer. Then one may prove that the derived space $d\ell_2$ is isomorphic to its own dual by \cite[Corollary 4]{Cab}. The second fact is the existence of a natural basis in such $d\ell_2$.
\begin{proposition}\label{lambda} Let $\mathbf X$ be a pseudolattice pair that is an $\ell_{\infty}$-module and denote by $(e_j)_{j=1}^{\infty}$ be the unconditional basis of $(B, B^*)_{\mathbf X, 1/2}=\ell_2$ and set $v_{2j-1}=(e_j,0), v_{2j}=(0,e_j)$ for $j\in \mathbb N$. Then
\begin{enumerate}
\item[$(i)$] The sequence $(v_j)_{j=1}^{\infty}$ is a basis for $d \ell_2$.
\item[$(ii)$] The sequence $(v_{2j})_{j=1}^{\infty}$ is unconditional.
\end{enumerate}
\end{proposition}
The first part is proved by adapting the proof of \cite[Theorem 4.10]{KaPe} while the second follows picking $a\in \{-1,1\}^{n}$ in (\ref{derivation}).
\section{The method of critical points}\label{sec4}
We give a simple way to compute derivations for a couple $(B,B^*)$. We shall only describe the derivations in the points, called informally ``critical points".
The technique applies both to the complex and the real method. 

For the real method, the CKMR-method produces for the pseudolattice pair $\lbrace\ell_p,\ell_p\rbrace$ a method equivalent to the $J$-method. It also agrees on the derivation operator \cite[Theorem 4.7]{CKMR}. Thus, for simplicity for us, we will work with the pseudolattice ${\lbrace \ell_2,\ell_2 \rbrace}$ instead of the $J$-method.

Hence, one of our applications deals with the derivation of the real interpolation formula $(T^2, (T^2)^*)_{\lbrace \ell_2,\ell_2 \rbrace,2}=\ell_2$ (cf. \cite{CoS}), we shall describe the method for this particular situation as an introductory example.
\begin{example} 
The critical points for $(T^2, (T^2)^*)_{\lbrace \ell_2, \ell_2 \rbrace ,2}$.
\end{example}
Let $\kappa(n)$ denote the equivalence constant of $(e_j)_{j=1}^n$ in $T^2$ with $\ell_2^n$. It is known that:
\begin{enumerate}
\item $\kappa(n)\to \infty$ as $n\to \infty$.
\item $$\lim_{n\to \infty}\frac{\kappa(n)}{\log_m(n)}=0,\;\;\;\;\textit{for}\;\;m=1,2,...$$
\end{enumerate}
The first claim follows since $T^2$ contains no copy of $\ell_2$. The second traces back to a close inspection of \cite[Proposition IV.b.3]{CS}. We are now ready to prove:
\begin{proposition}\label{extremal}
Let $\Omega$ be a derivation induced by $(T^2, (T^2)^*)_{\lbrace \ell_2, \ell_2 \rbrace ,2}$. There exists a sequence $(a_n)_{n=1}^{\infty}\subseteq \ell_2$ of non-null vectors with $\supp a_n \subseteq [1,n]$ such that 
\begin{equation}\label{quasilinear}
\Omega(a_n)=-2e^{-1/2}a_n\cdot \lfloor\log\kappa(n)\rfloor.
\end{equation}
\end{proposition}
\begin{proof} Let us recall that the norm of $\ell_2$ dominates the norm of $T^2$, so by definition of $\kappa(n)$ we have that 
\begin{equation*}
\|a\|_{T^2}\leq \|a\|_{\ell_2}\leq \kappa(n)\|a\|_{T^2}
\end{equation*}
for every $a\in[e_j]_{j=1}^n$, and by simple duality
\begin{equation*}
\kappa(n)^{-1}\|a\|_{(T^2)^*}\leq \|a\|_{\ell_2}\leq \|a\|_{(T^2)^*}.
\end{equation*}
Let us pick one ``critical" $a_n\in[e_j]_{j=1}^n$ satisfying $$\|a_n\|_{\ell_2} \sim  \kappa(n)\|a_n\|_{T^2}.$$ Let us give a bounded selection for this $a_n$ as follows: Pick $N=-2\lfloor \log\kappa(n)\rfloor$ and define
\begin{equation*}
b_j= \left\{ \begin{array}{lcc}
             e^{-N/2}a_n &,& j=N\\
             0 &,& j\neq N.\\
             \end{array}
   \right.
\end{equation*}
To see that it is a selection just write $$\sum_{j\in \mathbb Z}e^{\frac{j}{2}}b_j=e^{N/2}b_N=e^{N/2}e^{-N/2}a_n=a_n.$$
Let us show that the selection is bounded:
\begin{eqnarray*}\label{cota1}
\|b_N\|_{T^2}&=&\|e^{-N/2}a_n\|_{T^2}\\
&\leq & \| e^{\log\kappa(n)}a_n\|_{T^2}\\
&=&\kappa(n)\|a_n\|_{T^2}\\
&\sim& \|a_n\|_{\ell_2}.
\end{eqnarray*}
Similarly, we have
\begin{eqnarray*}
\|e^{N}b_N\|_{(T^2)^*}&=&\|e^{N/2}a_n\|_{(T^2)^*}\\
&=&\|e^{-\lfloor\log\kappa(n)\rfloor}a_n\|_{(T^2)^*}\\
&\leq & \| e^{1-\log\kappa(n)}a_n\|_{(T^2)^*}\\
&=&e\cdot \kappa(n)^{-1}\|a_n\|_{(T^2)^*}\\
&\leq& e\cdot \|a_n\|_{\ell_2}.
\end{eqnarray*}
Let us briefly observe that $$\|\{b_j\}_{j\in \mathbb Z}\|_{\ell_2(T^2)}=\|b_N\|_{T^2},$$  $$\|\{e^jb_j\}_{j\in \mathbb Z}\|_{\ell_2((T^2)^*)}=\|e^Nb_N\|_{(T^2)^*},$$ so that the selection is $e$-bounded for the pseudolattice $\mathcal X_0=\mathcal X_1=\ell_2$.

 Now, we may compute the derivation recalling the formal expression $$\Omega\left(\sum_{j\in \mathbb Z}e^{\frac{j}{2}}b_j\right)=e^{-1/2}\sum_{j\in \mathbb Z}je^{\frac{j}{2}}b_j$$ which gives us
$$\Omega(a_n)=e^{-1/2}Ne^{N/2}b_N=e^{-1/2}Na_n.$$
\end{proof}
\begin{corollary}\label{bounded}
$d(T^2, (T^2)^*)_{\lbrace \ell_2, \ell_2 \rbrace ,2}$ is not a Hilbert space.
\end{corollary}
\begin{proof}
If the derived space is isomorphic to $\ell_2$ then it has type $2$, and thus also $[v_{2j}]_{j=1}^{\infty}$ has type 2. But since $\{v_{2j}\}_{j=1}^{\infty}$ is an unconditional basis, it follows that the derivation operator $\Omega$ of (\ref{quasilinear}) is bounded (see \cite[Proposition 1]{JS2} for a detailed argument). But then, since $\log \kappa(n)$ is bounded, it follows that $\kappa(n)$ and $\kappa(n)^{-1}$ are bounded. In particular, $T^2$ and $(T^2)^*$ are isomorphic to $\ell_2$ which is absurd.
\end{proof}
\begin{example}The critical points working for the complex method.
\end{example}
Let us write the abstract version of this simple idea explicitly for the complex method. We shall assume some acquaintance of the reader with the complex method as introduced by Kalton and Montgomery-Smith \cite{KM}.
\begin{proposition}\label{propocomplex}
Let $(B,B^*)$ be a couple of sequence spaces. Assume that $(B,B^*)_{1/2}=\ell_2$ and let $\Omega$ be the corresponding derivation. Then, for each finite set $F\subseteq \mathbb N$, there is a nonzero $b_*\in [e_j]_{j\in F}$ such that \begin{equation}\label{centralizadorcom}
\Omega(b_*)=-2b_*\cdot \log \kappa(F),
\end{equation}
where 
$$\kappa(F)=\sup \{  \|b\|_{\ell_2}: \|b\|_B=1, \supp b\subseteq F \}.$$
\end{proposition}
\begin{proof}
By the very definition we have 
\begin{equation*}
\|b\|_{\ell_2}\leq \kappa(F)\|b\|_{B}
\end{equation*}
for every $b\in[e_j]_{j\in F}$, and by simple duality
\begin{equation*}
\kappa(F)^{-1}\|b\|_{B^*}\leq \|b\|_{\ell_2}.
\end{equation*}
Let us pick one critical $b_*\in[e_j]_{j\in F}$ satisfying $$\|b_*\|_{\ell_2} \sim  \kappa(F)\|b_*\|_{B},$$ that exists by compactness since $F$ is finite. Let us give a selection for this $b_*$ as follows: $$S_{b_*}(z)=e^{-2(z-1/2)\log\kappa(F)}b_*\in \mathcal F_{\infty}(B, B^*),$$
where $\mathcal F_{\infty}$ denotes the space of Calder\'on, see \cite{KM}.
Let us show that it is bounded:
\begin{eqnarray*}\label{cota2}
\|S_{b_*}(0+it)\|_{B}&=&\|e^{\log\kappa(F)}b_*\|_{B}\\
&\sim& \|b_*\|_{\ell_2}.
\end{eqnarray*}
Similarly, we have:
\begin{eqnarray*}
\|S_{b_*}(1+it)\|_{B^*}&=&\|e^{-\log\kappa(F)}b_*\|_{B^*}\\
&\leq & \|b_*\|_{\ell_2}.
\end{eqnarray*}
Therefore $\Omega(b_*)=\delta_{1/2}'S_{b_*}(z)$ and an easy evaluation finishes.
\end{proof}
An important situation where this can be used is the computation of the derivation corresponding to $Z(T^2)=d(T^2, (T^2)^*)_{1/2}$. A similar calculus as before shows that $$\Omega_{Z(T^2)}(a_n)=-2a_n\cdot \log \kappa(n).$$
This suggests that the spaces $Z(T^2)$ and $d(T^2, (T^2)^*)_{\lbrace \ell_2, \ell_2 \rbrace ,2}$ could be isomorphic but this is something we cannot prove. For the Kalton-Peck map, one may of course go further.  Consider the space $(c_0,\ell_{1})_{1/2}$ and pick $$a_n=\sum_{j=1}^n \frac{1}{\sqrt{n}}e_j$$ that is critical and for which the equivalence constant is $\sqrt{n}$. This gives us back, of course, the Kalton-Peck map on $a_n$:
$$\Omega(a_n)=-2a_n\cdot \log \sqrt{n}.$$
We would like to remark that the norm of $T^2$ is given implicitly, so perhaps one cannot expect to have an explicit formula of $\Omega_{Z(T^2)}$ for all the points as in the Kalton-Peck map.
\subsection{A geometric interpretation}\label{geometric}
Let us give a geometric interpretation of these results using ideas of Semmes \cite[Page 160]{Se}. Let $\mathcal N(\mathbb C^n)$ be the set of all norms over $\mathbb C^n$. Complex interpolation theory provides a sort of structure over this set that make it looks like a Riemannian manifold in some ways. Given points $\|\cdot\|_0,\|\cdot\|_1\in \mathcal N(\mathbb C^n)$, we define a curve joining them by the formula $$\|\cdot\|_{\theta}=(\|\cdot\|_0,\|\cdot\|_1)_{\theta},\;\;\;\;\;\;\theta\in[0,1].$$ 
Semmes reads these ``Calder\'on curves" as geodesics that minimize length. This length is given in a natural way for a couple $M,N\in \mathcal N(\mathbb C^n)$. First, we introduce a gap $$\delta(M,N)=\log \sup_{0\neq v\in \mathbb C^n} \frac{M(v)}{N(v)},$$
and then the distance is given as $$d_{\mathcal C}(M,N):=\sup \{\delta(M,N), \delta(N,M)\}.$$
If we restrict ourselves to vectors with support in $[1,n]$, the derivation introduced in (\ref{centralizadorcom}) plays the role of a tangent vector to the Calder\'on curve joining $B_n=[e_j]_{j=1}^n$ and $B^*_n=[e_j^*]_{j=1}^n$, where $(e_j)_{j=1}^{\infty}$ and $(e_j^*)_{j=1}^{\infty}$ denotes the corresponding basis in $B$ and $B^*$. Assume for simplicity that the norm of $B^*$ dominates the norm of $B$. Therefore, since $\ell_2$ is an intermediate space, we have  $\|\cdot\|_B\leq  \|\cdot\|_{\ell_2} \leq \|\cdot\|_{B^*}$. In particular, $$\|\cdot\|_B\leq  \|\cdot\|_{\ell_2} \leq \kappa([1,n])\|\cdot\|_{B}$$ for vectors with support in $[1,n]$. Thus,
$$d_{\mathcal C}(B_n,\ell_2^n)=\log\kappa([1,n]).$$
So, we may rewrite this tangent vector (\ref{centralizadorcom}) as:
$$\Omega(b_n)=-2b_n \cdot d_{\mathcal C}(B_n,\ell_2^n),$$
where $b_n\neq 0$ is the critical point with support in $[1,n]$ guaranteed by Proposition \ref{propocomplex}. So the derivation depends only of the ``Calder\'on distances" between the norms of $B$ and $\ell_2$. An analogous formula involving $d_{\mathcal C}(B_n^*,\ell_2^n)$ holds true, of course. One just need to introduce, say $\kappa^*(F)$, in a similar vein. This is done in the next section where we describe the whole picture in an abstract way.
\section{A universal formula for derivation operators}\label{sec5}
We introduce the formula for derivations sketched in the previous section in full generality and close this section with a version of Kalton uniqueness theorem for a large class of pseudolattice pairs $\mathbf X$: those for which the norm of a sequence $\{b_n\}_{n\in \mathbb Z}$ supported in one single element, say $b_{n_0}\in B\cap B^*$, satisfies
\begin{equation}\label{bo}
 \left\{ \begin{array}{lcc}
             \|\{b_n\}_{n\in \mathbb Z}\|_{\mathcal X_0(B)}\sim \|b_{n_0}\|_{B},\\
              \|\{b_n\}_{n\in \mathbb Z}\|_{\mathcal X_1(B^*)}\sim \|b_{n_0}\|_{B^*}.\\
             \end{array}
   \right.
\end{equation}
The pseudolattices $FC$ and $\ell_p$ do trivially satisfy the condition above, so the results of this section generalize to those of the previous section for the complex and the real method, respectively.
Given a finite subset $F\subseteq \mathbb N$, we define
$$\kappa^*(F)=\sup \{  \|b\|_{\ell_2}: \|b\|_{B^*}=1, \supp b\subseteq F \}.$$ We are ready to prove now the formula in full generality.
\begin{theorem}
Let $(B,B^*)$ be a couple of sequence spaces with basis $(e_j)_{j=1}^{\infty}$ and let $\mathbf X$ be a pseudolattice pair satisfying property (\ref{bo}). Assume that $(B,B^*)_{\mathbf X, 1/2}=\ell_2$ and let $\Omega$ be the corresponding derivation operator. Then, for each finite set $F\subseteq \mathbb N$, there are  nonzero $b, b^*\in [e_j]_{j\in F}$ such that \begin{equation}\label{canonical}
\Omega(b)=-2e^{-1/2}b\cdot\lfloor\log \kappa(F)\rfloor,
\end{equation}
\begin{equation}\label{canonical2}
\Omega(b^*)=+2e^{-1/2}b^*\lfloor\log \kappa^*(F)\rfloor.
\end{equation}
\end{theorem}
\begin{proof}
The proof is similar to the given in Proposition \ref{extremal}. Let us fix a finite set $F$ for the rest of the proof. Observe that by definition we have that $\|b\|_{\ell_2}\leq \kappa(F)\|b\|_B$ and $\|b\|_{\ell_2}\leq \kappa^*(F)\|b\|_{B^*}$ for all $b$ with $\supp b\subseteq F$. Then, by simple duality, we also have for all such $b$ that $$\kappa(F)^{-1}\|b\|_{B^*}\leq \|b\|_{\ell_2}$$
and
$$\kappa^*(F)^{-1}\|b\|_{B}\leq \|b\|_{\ell_2}.$$
Pick $b^*$ such that $$\|b^*\|_{\ell_2}\sim \kappa^*(F)\|b^*\|_{B^*}$$ and define for $N=2\lfloor\log \kappa^*(F)\rfloor$
\begin{equation*}
b_j= \left\{ \begin{array}{lcc}
            e^{-N/2}b^* &,& j=N\\
             0 &,& j\neq N.\\
             \end{array}
   \right.
\end{equation*}
A similar computation as given in Proposition \ref{extremal} shows that it is a bounded selection from where the second claim follows. For the first claim, pick $b$ such that $$\|b\|_{\ell_2}\sim \kappa(F)\|b\|_{B},$$
and argue exactly as in Proposition \ref{extremal}.
\end{proof}

We close with our version of Kalton uniqueness theorem for a large class of methods. We would like to recall that Kalton proved that if the derivation operator $\Omega$ corresponding to a couple $(B_0,B_1)$ satisfies that, for some constant $C$ and certain linear map $L$, $$\|\Omega b-Lb \| \leq C \|b\|,$$ then $B_0$ and $B_1$ have equivalent norms. We prove a much weaker version of this result where the derivation operator is the corresponding to the Hilbert space for a couple $(B,B^*)$ and $L=0$.
\begin{theorem}Let $(B,B^*)$ be a couple of sequence spaces. Let $\mathbf X$ be a pseudolattice pair satisfying property (\ref{bo}). Assume that $(B,B^*)_{\mathbf X, 1/2}=\ell_2$ and let $\Omega$ be the corresponding derivation operator. If $\Omega$ is bounded then the norms of $B$ and $B^*$ are equivalent to the norm of $\ell_2$.
\end{theorem} 
\begin{proof}
If $\Omega$ is bounded, then so does $\kappa(F)$ and $\kappa^*(F)^{-1}$ (and also $\kappa^*(F)$ and $\kappa(F)^{-1}$) uniformly for all the finite sets $F\subseteq \mathbb N$, which gives that the norm of $\ell_2$ is equivalent to $B$ (and also to $B^*$).
\end{proof}
\section{Applications}
In this section we give two examples of twisted Hilbert spaces that now also arise by real methods and we close constructing derivation operators whose norm on critical points grow as slowly as we want.
\subsection{On real interpolation methods: a word}\label{word}
There are many interpolation methods in the literature considered as ``real": the $K/J$-methods the most noticeable but also the $E$ method \cite[Chapter 7]{BeLo}, its dual counterpart the $F$ method, the Lions-Peetre \textit{espaces de moyennes} or some choices described in \cite{CKMR}. All these methods produce the same interpolation space, up to equivalence of norms. The delicate point here is that they do not produce necessarily the same derived space. For instance, the $K$-method produces, for the usual pair $(L_1,L_{\infty})$, the derivation operator $$\Omega_{K}(f)=f|\log r_f|,$$
for a normalized $f\in L_p$ and where $r_f$ denotes the rank function, see for example \cite{JRW}. However, for such normalized $f\in L_p$, the $E$-method gives the derivation map $$\Omega_E(f)=f\log|f|,$$
which is the Kalton-Peck map. A proof may be found at \cite{JRW}. So, the precise choice of ``real" method is crucial in the context of derivation maps; even in the setting of natural couples as $(\ell_{p_0},\ell_{p_1})$. Indeed, it is very well known that for the $J$-method we always have $\Omega_K=-\Omega_J$ (\cite{CJM}), so we do find for the couple $(\ell_{p_0},\ell_{p_1})$ that $$\Omega_J(x)=-x|\log r_x|,$$
for a normalized $x\in \ell_p$ and where $r_x$ denotes again the rank function. The CKMR-method produces for the pseudolattice pair $\lbrace\ell_p,\ell_p\rbrace$ a method equivalent to the $J$-method. It also agrees on the derivation operator \cite[Theorem 4.7]{CKMR}. Thus, for the couple $(\ell_{p_0},\ell_{p_1})$, with $p_0\neq p_1$, we have 
\begin{equation}\label{toma}\Omega_{\lbrace\ell_p,\ell_p\rbrace,\theta}(x)=-x|\log r_x|,
\end{equation}
for such normalized $x\in \ell_p$ and where $(1-\theta)p_0^{-1}+\theta p_1^{-1}=1/p$. The CKMR-method for the pseudolattice $\lbrace\ell_{p_0},\ell_{p_1}\rbrace$ and our fixed $\theta$ as above is also a real method equivalent to the $J$-method of parameters $(\theta,p)$. This is argued in the last paragraph before \cite[Section 3]{CKMR} and means that both methods produce isomorphic Banach spaces. However, Theorem \ref{Lions} (in the next section) will show that this real CKMR-method with the pseudolattice $\lbrace\ell_{p_0},\ell_{p_1}\rbrace$ produces now
\begin{equation}\label{hawk}\Omega_{\lbrace \ell_{p_0},\ell_{p_1}\rbrace, \theta}(x)=c\cdot x \lfloor\log x \rfloor,
\end{equation}
for such normalized $x\in \ell_p$ and where $c$ is a numerical constant depending only on $\theta, p_0, p_1$.

%
%

To put it tersely, the two real CKMR-methods given by $(\lbrace\ell_{p},\ell_{p}\rbrace, \theta)$ and $(\lbrace\ell_{p_0},\ell_{p_1}\rbrace, \theta)$, where $1/p=(1-\theta)p_0^{-1}+\theta p_1^{-1}$ and $p_0\neq p_1$, produce isomorphic Banach spaces for any interpolation couple but different derivation maps: (\ref{toma}) and (\ref{hawk}), respectively.

\subsection{A ``real" Kalton-Peck space}\label{kaltonpeque}
That the Kalton-Peck space arises also from the real method is implicit in the existing literature as quoted above. We have tracked it to the pioneering paper of Lions and Peetre \cite{LiPe}. The reason for this is clear: Lions and Peetre computed the extremal in the proof of \cite[Theorem (I.I), Chapitre VII]{LiPe} in an easy and elegant way. The equivalence between the discrete version of the Lions-Peetre \textit{espaces de moyennes} and the approach of Cwikel \textit{et al.} is quite natural. One just need to pass the extremal for this natural equivalence. In practice, this means to multiply by the weight $e^{-n\theta}$. Once this is done, we just need to make a simple substitution in the derivation formula. 
\begin{theorem}(Lions-Peetre)\label{Lions}
The derivation operator associated to $$(\ell_{p_0}, \ell_{p_1})_{\lbrace\ell_{p_0},\ell_{p_1}\rbrace,\theta}=\ell_p$$ is given for $a=(a_m)_{m\in \mathbb Z}\in \ell_p$ as 
\begin{equation}\label{kpreal}
\Omega(a)=e^{-\theta}\sum_{m\in \mathbb Z}  -\left\lfloor\left(\frac{p}{p_0}-\frac{p}{p_1}\right)\log \frac{|a_m|}{\|a\|_p} \right\rfloor a_me_m ,
\end{equation}
where $\frac{1}{p}=\frac{1-\theta}{p_0}+\frac{\theta}{p_1}$, $0<\theta<1$ and $1\leq p_0,p_1<\infty$.
\end{theorem}
\begin{proof}
The formula $\ell_p=(\ell_{p_0},\ell_{p_1})_{\lbrace\ell_{p_0},\ell_{p_1}\rbrace,\theta}$ is established in \cite[Theorem (I.I), Chapitre VII]{LiPe}. Therefore, for a fixed $C>1$, we have that
\begin{equation}\label{formula}
\Omega(a)=\sum_{n\in \mathbb Z} n e^{\theta(n-1)}b_n,
\end{equation}
where \begin{equation}\label{section}
\delta_{\theta}(\{b_n\}_{n\in \mathbb Z})=a
\end{equation} with
 \begin{equation}\label{constant}
 \|\{b_n\}_{n\in \mathbb Z}\|\leq C \|a\|_{\ell_p}.
 \end{equation}
So one just need to compute $\{b_n\}_{n\in \mathbb Z}$ in terms of $a=\{a_m\}_{m\in \mathbb Z}$. To this end let us find $\lambda$ such that
\begin{equation*}
\left\{ \begin{array}{lcc}
            p_0(1+\lambda \theta)=p, \\
             p_1(1-\lambda(1-\theta))=p,\\
             \end{array}
   \right.
\end{equation*}
 exactly as in the proof of \cite[Theorem (I.I), Chapitre VII]{LiPe}. Indeed, one may take 
\begin{equation*}\label{lam}
\lambda=\frac{p}{p_0}-\frac{p}{p_1}.
\end{equation*}
There is no loss of generality to assume $\|a\|=1$ and then extend by homogeneity. Lions and Peetre defined the corresponding $\{b_n\}_{n\in \mathbb Z}$ for a given $a\in \ell_p$ with $\|a\|=1$ as follows
\begin{equation}\label{b}
b_n(m)= \left\{ \begin{array}{lcc}
             e^{-n\theta}a_m &,& n=-\lfloor\lambda \log |a_m|\rfloor\\
             \\ 0 &,& \textit{otherwise.}\\
             \end{array}
   \right.
\end{equation}
It is trivial to check that condition (\ref{section}) holds. Let us check that condition (\ref{constant}) is satisfied:
\begin{eqnarray*}
\|\{b_n\}_{n\in \mathbb Z}\|_{\ell_{p_0}(\ell_{p_0})}^{p_0}&=&\sum_{n} \|b_n\|^{p_0}_{\ell_{p_0}}\\
&=& \sum_{n=-\lfloor\lambda \log |a_m|\rfloor} \left |e^{-n\theta}a_m \right|^{p_0}\\
& \leq & \sum_m e^{\theta(\lambda \log |a_m|)p_0} |a_m|^{p_0}\\
&=& \sum_m  |a_m|^{p_0+\lambda \theta p_0}\\
&=&\sum_m |a_m|^p=1.
\end{eqnarray*}
In a similar vein one may prove that $$\|\{e^nb_n\}_{n\in \mathbb Z}\|_{\ell_{p_1}(\ell_{p_1})}\leq 1.$$ Thus, the condition  (\ref{constant}) is satisfied with $C=1$. A simple evaluation in the formula (\ref{formula}) for the values of $\{b_n\}_{n\in \mathbb Z}$ given in (\ref{b}) gives the desired result; observe that the non void terms are $n=-\lfloor\lambda \log |a_m|\rfloor$.
\end{proof}
\subsection{A ``real" weak Hilbert space}\label{sec3}
Our aim is to show that $d(T^2, (T^2)^*)_{\lbrace\ell_{2},\ell_{2}\rbrace, 2}$ is a weak Hilbert space. As in \cite{JS}, we need to estimate first the type $p$ constant of the derived space and this requires a technical result on random sums for derivations. The following lemma is well known.
\begin{lemma}\label{interpola}
Let $1<p\leq 2$ and $p\leq q<\infty$. If $B=(B_0,B_1)$ is a compatible pair, then, for every $m\geq 1$ we have
 $$a_{m,p}(B_{\lbrace\ell_{q},\ell_{q}\rbrace,\theta})\lesssim a_{m,p}(B_0)^{1-\theta}a_{m,p}(B_1)^{\theta}.$$
\end{lemma}
\begin{proof}
Since the shift operator is an isometry acting on $\ell_q(B_j)$ with $j=0,1$, then we have $$\|\cdot\|_{\theta,q}\leq e\cdot \|\cdot\|_{\ell_q(B_0)}^{1-\theta}\|\cdot\|_{\ell_q(B_1)}^{\theta},$$
for $1\leq q \leq \infty$, see \cite[Lemma 2.13]{CKMR}. Now average over all $\varepsilon\in \{-1,1\}^m$ and use Holder's inequality. To conclude use that $a_{m,p}(\ell_q(B_k))\lesssim a_{m,p}(B_j)$ for $j=0,1$, see \cite[Theorem 11.12(a)]{DJT}.
\end{proof}
We are ready to introduce the technical result on random sums.
\begin{lemma}\label{randoma}
Let $1<p\leq 2$ and $p\leq q<\infty$. Let $B=(B_0,B_1)$ be a compatible pair and $\Omega$ a derivation induced by $B_{\lbrace\ell_{q},\ell_{q}\rbrace,\theta}$ with $0<\theta<1$. Then, for every $b_1,...,b_m\in B_{\lbrace\ell_{q},\ell_{q}\rbrace,\theta}$, we have
that 
\begin{equation}\label{basic}
\mathbb E \left\| \Omega\left(\sum_{j=1}^m \varepsilon_j b_j\right)-\sum_{j=1}^m \varepsilon_j\Omega(b_j)-e^{-1}\log \frac{a_{m,p}(B_0)}{a_{m,p}(B_1)}\sum_{j=1}^m \varepsilon_j b_j\right\|
\end{equation} 
is upper bounded, up to a constant $C(p,q)$ depending only on $p,q$, by $$ a_{m,p}(B_{0})^{1-\theta}a_{m,p}(B_{1})^{\theta} \left(\sum_{j=1}^{m}\|b_j\|^p\right)^{1/p}.$$
\end{lemma}
\begin{proof} Recall that the pseudolattice pair $\mathbf X=\lbrace\ell_q,\ell_q\rbrace$ has it corresponding ``Calder\'on space", namely $\mathcal J(\mathbf X, B)$. Fix $S$ an homogeneous $2$-bounded selector for the quotient map $\delta_{\theta}$ onto $\mathcal J(\mathbf X, B)$ and let $\Omega:=\delta'_{\theta}S$ be the corresponding derivation. Given $b_j$, let us write informally
\begin{equation}\label{sel}
S(b_j)=\sum_{n\in \mathbb Z}b_{j,n}z^n, \;\;\;\;j=1,...,m.
\end{equation}
To simplify the forthcoming computations let us denote simply $a_k=a_{m,p}(B_k)$ for $k=0,1$. For each $\varepsilon\in \lbrace -1,1 \rbrace^m$, we define the map 
$$F_{\varepsilon}(z)=a_0^{-1}\sum_{j=1}^m\varepsilon_j\sum_{n\in \mathbb Z} b_{j,n}z^{n+\log a_0-\log a_1}.$$
Formally speaking we should deal with the map
$$\lfloor F_{\varepsilon}\rfloor(z)=e^{-\lfloor\log a_0\rfloor}\sum_{j=1}^m\varepsilon_j\sum_{n\in \mathbb Z} b_{j,n}z^{n+\lfloor\log a_0\rfloor-\lfloor\log a_1\rfloor}.$$
More precisely, this map comes from the following. Given $b_j$ for $j\leq m$, we have formally $S(b_j)=\{b_{j,n}\}_{n\in \mathbb Z}\in \mathcal J(\mathbf X,B)$. Now, let us define for each $n\in \mathbb Z$: $$\widetilde{b}_{j,n}=b_{j,n-\lfloor\log a_0\rfloor+\lfloor\log a_1\rfloor}.$$
Then we have that 
\begin{eqnarray*}
\left\{e^{kn}\widetilde{b}_{j,n}\right\}_{n\in \mathbb Z}&=&\left\{e^{kn}b_{j,n-\lfloor\log a_0\rfloor+\lfloor\log a_1\rfloor}\right\}_{n\in \mathbb Z}\\
&=&\left\{e^{k(n+\lfloor\log a_0\rfloor-\lfloor\log a_1\rfloor)}b_{j,n}\right\}_{n\in \mathbb Z}\in \ell_q(B_k),\;\;\;k=0,1.\\
\end{eqnarray*}
In particular, 
\begin{eqnarray*}
&\left\{e^{-\lfloor\log a_0\rfloor}\sum_{j=1}^m \varepsilon_j \widetilde{b}_{j,n}\right\}_{n\in \mathbb Z}\in \mathcal J(\mathbf X,B).\;\;\;\\
\end{eqnarray*}
The map $\lfloor F_{\varepsilon}\rfloor$ corresponds to the element above but for simplicity we will use the map $F_{\varepsilon}$.
%
%
Observe now that we have
\begin{eqnarray*}
a_0^{1-\theta}a_1^{\theta}F_{\varepsilon}'(e^{\theta})&=& a_0^{1-\theta}a_1^{\theta} \left[\sum_{j=1}^m \varepsilon_j \sum_{n\in \mathbb Z} b_{j,n}a_0^{-1}(n+\log a_0-\log a_1)e^{\theta(n-1+\log a_0-\log a_1)} \right]\\
&=& a_0^{1-\theta}a_1^{\theta}\left[ \sum_{j=1}^m \varepsilon_j \sum_{n\in \mathbb Z} b_{j,n}a_0^{-1}(n+\log a_0-\log a_1)e^{\theta(n-1)}a_0^{\theta}a_1^{-\theta} \right]\\
&=&\sum_{j=1}^m \varepsilon_j\sum_{n\in \mathbb Z} b_{j,n}ne^{\theta(n-1)}+e^{-1}\log\frac{a_0}{a_1}\sum_{j=1}^m \varepsilon_j \sum_{n\in \mathbb Z}  b_{j,n}e^{\theta n}\\
&=& \sum_{j=1}^m \varepsilon_j \Omega(b_j)+e^{-1}\log\frac{a_0}{a_1}\sum_{j=1}^m \varepsilon_j b_j,
\end{eqnarray*}
which corresponds to the ``linear term" in the inequality (\ref{basic}) of the theorem. Indeed, we should have:
$$e^{(1-\theta)\lfloor\log a_0\rfloor}e^{\theta\lfloor\log a_1\rfloor}\lfloor F_{\varepsilon}\rfloor'(e^{\theta})= \sum_{j=1}^m \varepsilon_j \Omega(b_j)+e^{-1}(\lfloor\log a_0\rfloor-\lfloor\log a_1\rfloor)\sum_{j=1}^m \varepsilon_j b_j.$$
 Once this is done, we may reduce the argument to the ``cancellation principle" as in \cite[Lemma 1]{JS}. 
\begin{eqnarray*}
\Lambda_{\varepsilon}&=&\left \| \Omega\left(\sum_{j=1}^m \varepsilon_j b_j\right)-\sum_{j=1}^m \varepsilon_j\Omega(b_j)-e^{-1}\log \frac{a_{m,p}(B_0)}{a_{m,p}(B_1)}\sum_{j=1}^m \varepsilon_j b_j   \right \|_{\theta, q}\\
&=& \left\| \delta_{\theta}'S\left(\sum_{j=1}^m \varepsilon_j b_j \right)- \delta_{\theta}'\left(a_0^{1-\theta}a_1^{\theta}F_{\varepsilon}(z)  \right) \right\|_{\theta, q} \\
&=& \left\| \delta_{\theta}' \left( S \left(\sum_{j=1}^m \varepsilon_j b_j \right)-a_0^{1-\theta}a_1^{\theta}F_{\varepsilon}(z) \right) \right\|_{\theta, q}  \\
\end{eqnarray*}
A critical observation is that $S \left(\sum_{j=1}^m \varepsilon_j b_j \right)-a_0^{1-\theta}a_1^{\theta}F_{\varepsilon}(z)\in \Ker \delta_{\theta}$. Recall that the shift operator is an isometry on $\ell_q(B_k)$ for $k=0,1$. Therefore, we may apply \cite[Lemma 3.6]{CKMR} to conclude that our pseudolattice \textit{admits differentation} (\cite[Definition 3.4]{CKMR}) and thus apply \cite[Lemma 3.11]{CKMR}. This last result combined with the same argument given in \cite[Lemma 1, line 11]{JS}, gives immediately that there is a numerical constant $C=C(\theta)$ so that 
$$\Lambda_{\varepsilon}\leq C\left(\left\| S \left(\sum_{j=1}^m \varepsilon_j b_j \right)\right \| + a_0^{1-\theta}a_1^{\theta}\|F_{\varepsilon}(z)\|\right),$$
so we are led to bound them separately. For the first term, just write
$$\mathbb E \left\| S \left(\sum_{j=1}^m \varepsilon_j b_j \right)\right \|\leq 2 \mathbb E \left\| \sum_{j=1}^m \varepsilon_j b_j\right \|\leq 2a_{m,p}(B_{\lbrace \ell_q,\ell_q \rbrace,\theta})\left(\sum_{j=1}^{m}\|b_j\|_{\lbrace \ell_q,\ell_q \rbrace, \theta}^p\right)^{1/p},$$
and use Lemma \ref{interpola}. For the second term, first write informally
\begin{eqnarray*}
\mathbb E \left \| F_{\varepsilon}(z)  \right\|&=& \mathbb E \left( \max_{k=0,1}\left \| F_{\varepsilon}(e^k) \right\|_{\ell_q(B_k)}\right)\\
&\leq& \mathbb E\left \| F_{\varepsilon}(e^0) \right\|_{\ell_q(B_0)}+\mathbb E\left \| F_{\varepsilon}(e^1) \right\|_{\ell_q(B_1)}.
\end{eqnarray*}
To bound each term of the sum above let us recall something well known. The space $L_q(\mu, B_k)$ has type $p\leq 2$ if (and only if) $B_k$ has type $p$, whenever $p\leq q$, \cite[Theorem 11.12(a)]{DJT}. Indeed the proof of this fact shows that there is $c=c(p,q)$ so that for every $m\geq 1$ $$a_{m,p}(L_q(\mu,B_k))\leq c(p,q)a_{m,p}(B_k).$$
Put this card in the pocket for a second and let us bound the terms of the sum
\begin{eqnarray*}
\mathbb E \left \| F_{\varepsilon}(e^k)  \right\|&=& a_k^{-1}\mathbb E\left \| \sum_{j=1}^m \varepsilon_j \sum_{n\in \mathbb Z} e^{kn}b_{j,n}  \right \|_{\ell_q(B_k)}\\
&\leq&\frac{a_{m,p}(\ell_q(B_k))}{a_k}\left( \sum_{j=1}^m  \left\|\sum_{n\in \mathbb Z} e^{kn}b_{j,n}  \right\|^p_{\ell_q(B_k)}  \right)^{1/p}\\
&\leq& c(p,q)\left( \sum_{j=1}^m  \left\|\sum_{n\in \mathbb Z} e^{kn}b_{j,n}  \right\|^p_{\ell_q(B_k)}  \right)^{1/p}\\
&\leq &2c(p,q) \left( \sum_{j=1}^m \|b_j\|^p_{\theta,q} \right)^{1/p},
\end{eqnarray*}
where we have used the card in the second inequality. And the proof of the main argument is complete. Formally, we have ended up with a $e^{-1}(\lfloor\log a_0\rfloor-\lfloor\log a_1\rfloor)$ factor in the linear term (\ref{basic}). This can be replaced by $e^{-1}\frac{\log a_0}{\log a_1}$ since 
\begin{equation*}
\left|e^{-1}(\lfloor\log a_0\rfloor-\lfloor\log a_1\rfloor)-e^{-1}\frac{\log a_0}{\log a_1}\right|\leq 2e^{-1}.
\end{equation*}
\end{proof}
Now, we may estimate the type $p$ constants of the derived space by the real  $(\lbrace\ell_q,\ell_q\rbrace,\theta)$-method.
\begin{proposition}\label{average}
Let $1<p\leq 2$ and $p\leq q<\infty$. Let $B=(B_0,B_1)$ be a compatible pair and $\Omega$ a derivation induced by $B_{\lbrace\ell_q,\ell_q\rbrace,\theta}$ with $0<\theta<1$. Then there is $C=C(\theta,p,q)$ such that for every $m\geq 1$
\begin{equation}
a_{m,p}(dB_{\lbrace\ell_q,\ell_q\rbrace,\theta})\leq C \cdot a_{m,p}(B_{\lbrace\ell_q,\ell_q\rbrace,\theta})\left(a_{m,p}(B_{0})^{1-\theta}a_{m,p}(B_{1})^{\theta} + \left|\log \frac{a_{m,p}(B_0)}{a_{m,p}(B_1)}\right|\right).
\end{equation}
\end{proposition}
\begin{proof}
It is exactly as \cite[Proposition 3]{JS} using Lemma \ref{randoma}. We reproduce the argument. Pick vectors $x_j=(a_j,b_j)$ in $dB_{\lbrace\ell_q,\ell_q\rbrace,\theta}$ for $j\leq m$ and let us denote by simplicity
$$\Delta_{\varepsilon}=\left \|\Omega \left(\sum_{j=1}^m \varepsilon_j b_j\right)- \sum_{j=1}^m \varepsilon_j\Omega(b_j)\right\|,\;\; \varepsilon\in \{-1,1\}^m.$$ Then 
\begin{eqnarray*}
\mathbb E\left \| \sum_{j=1}^m \varepsilon_j x_j \right \| &\leq & \mathbb E \left  \| \sum_{j=1}^m \varepsilon_j(a_j - \Omega(b_j))\right\| + \mathbb E\left \|  \sum_{j=1}^m \varepsilon_j b_j \right \|+ \mathbb E  \Delta_{\varepsilon}\\
&\leq & 2a_{m,p}(B_{\lbrace\ell_q,\ell_q\rbrace,\theta})\left(  \sum_{j=1}^m \|x_j\|^p\right) ^{1/p}+\mathbb E  \Delta_{\varepsilon}.
\end{eqnarray*}
A quick computation finishes the proof since we have  by Lemma \ref{randoma} $$\mathbb E  \Delta_{\varepsilon} \leq \gamma \cdot\left(\sum_{j=1}^m\|b_j\|^p\right)^{1/p}+\mathbb E\left\|\log \frac{a_{m,p}(B_0)}{a_{m,p}(B_1)} \sum_{j=1}^m \varepsilon_j b_j\right\|,$$
where $$\gamma \lesssim a_{m,p}(B_0)^{1-\theta}a_{m,p}(B_1)^{\theta}.$$
\end{proof}
%
 Let us recall that $d(T^2, (T^2)^*)_{\lbrace\ell_2,\ell_2\rbrace,1/2}$ is not a Hilbert by Corollary \ref{bounded}, so let us show now that 
\begin{proposition}
The derived space $d(T^2, (T^2)^*)_{\lbrace\ell_2,\ell_2\rbrace,\frac{1}{2}}$ is a weak Hilbert space.
\end{proposition}
\begin{proof}
The proof is similar to the one given for $Z(T^2)$ in \cite{JS}. Let $(e_j)$ and $(e_j^*)$ be the bases in $T^2$ and $(T^2)^*$ respectively and let us write $E_n=[e_j]_{j=n}^{\infty}$ and $E_n^*=[e_j^*]_{j=n}^{\infty}$. Once this is done, it is not hard to see that $$V_n=d(E_n,E_n^*)_{\lbrace\ell_2,\ell_2\rbrace,\frac{1}{2}}$$
where $V_n=[v_j]_{2n-1}^{\infty}$. The space $d(T^2, (T^2)^*)_{\lbrace\ell_2,\ell_2\rbrace,1/2}$ is $\lambda$-isomorphic to its own dual cf. \cite{Cab}. It is therefore plain to deduce that $V_n$ is also $\lambda$-isomorphic to its own dual for every $n\in \mathbb N$. Once we are equipped with this, the proof consists in three steps. For the first one, just replace \cite[Proposition 3]{JS} by Proposition \ref{average} with $p=q=\theta^{-1}=2$ using the previous indications while for the second and third steps just mimic the proof of \cite[Proposition 3]{JS}. 
\end{proof}

Let us denote $(T^2, (T^2)^*)_{\frac{1}{2},2;J}$ the interpolation space by the real $J$-method with parameters $(1/2, 2)$. Then, the result above is equivalent to say 

\begin{proposition}
The twisted Hilbert space $d(T^2, (T^2)^*)_{\frac{1}{2},2;J}$ is a weak Hilbert space.
\end{proposition}
\begin{proof}One just need to apply \cite[Theorem 4.7]{CKMR}, see also the comments of Subsection \ref{word}.
\end{proof}

Very recently, two new twisted Hilbert spaces have been introduced in \cite{JS3}, namely, $Z(\mathcal J)$ and $Z(\mathcal T_s^2)$. The first is asymptotically Hilbertian but not weak Hilbert while the second has the hereditary approximation property but it is not asymptotically Hilbertian. The definitions of both notions can be found at the first two paragraphs of \cite[Introduction]{JoSa}. It can be proved with the techniques developed here so far that both spaces do admit its ``real" counterpart by the $J$-method.

We would like to mention that throughout the paper the word ``real" makes reference to the method of interpolation and not to the scalar field. Indeed, the CKMR-method uses complex scalars so all the derived spaces appearing throughout the text are complex Banach spaces.
\subsection{Derivations \textit{\`a la} Johnson-Szankowski}
 Let us give one last application of our ideas. Until here, we have computed the maximum values of a given derivation. Now, we find a derivation operator once the maximum values are given. We need to recall results of Johnson and Szankowski (cf. \cite{JoSa}) where they construct, for a given sequence $1\leq \delta_n\longrightarrow \infty$, a Banach sequence space $X$ so that $d_n(X)\leq \delta_n$, where $$d_n(X)=\sup d(E,\ell_2^n),$$
the supremum runs over all the $n$-dimensional subspaces $E$ of $X$ and $d$ stands as usual for the Banach-Mazur distance. We may construct nontrivial derivation operators, even a centralizer, whose ``maximum values" on the first $n$-coordinates tends to infinity as slowly as we wish.
\begin{proposition} Let $1<\delta_n\to  \infty$. Then there is a nontrivial centralizer $\Omega:\ell_2 \longrightarrow \omega$ and a sequence $(a_n)_{n=1}^{\infty}\subseteq \ell_2$ of non-null vectors with $\supp a_n\subseteq [1,n]$ such that
$$\Omega(a_n)=-2a_n\cdot \log \kappa(n),$$
with 
$$\delta_n\geq \kappa(n) \longrightarrow \infty.$$
\end{proposition}
\begin{proof}
Use the Johnson-Szankowski space $X$ for the choice $(\delta_n)_{n=1}^{\infty}$, so that $(X,X^*)_{1/2}=\ell_2$ by \cite{CoS}. It is not hard to see that the corresponding derivation satisfies the claim of the proposition.
\end{proof}
The reader should perhaps compare this result with the Kalton-Peck derivation where $\kappa(n)=\sqrt{n}$ which gives the largest possible Calder\'on distance on $\mathcal N(\mathbb C^n)$ to the Hilbert norm.  In this sense, proposition above shows that there is no lower bound.


\begin{thebibliography}{99}

\bibitem{AK} F. Albiac and N. J. Kalton, \emph{Topics in Banach space theory}.
Graduate Texts in Mathematics 233. Springer-Verlag 2016.

\bibitem{BeLo}  J. Bergh,and J. L\"ofstr\"om, \emph{Interpolation spaces. An introduction}. Grundlehren der Mathematischen Wissenschaften, No. 223. Springer-Verlag, Berlin-New York, 1976. 


\bibitem{Cab} F. Cabello S\'anchez, {\it Nonlinear centralizers in homology}, Math. Ann. 358 (2014), no. 3-4, 779--798.


\bibitem {CS} P. Casazza and T.J. Shura, \emph{Tsirelson's space. With an appendix by J. Baker, O. Slotterbeck and R. Aron}. Lecture Notes in Mathematics, 1363. Springer-Verlag, Berlin, 1989.



\bibitem{CKMR} M. Cwikel, N. Kalton, M. Milman and R. Rochberg, {\it A unified theory of commutator estimates for a class of interpolation methods}, Adv. Math. 169 (2002), no. 2, 241--312. 

\bibitem{CJMR} M. Cwikel, B. Jawerth, M. Milman and R. Rochberg, {\it Differential estimates and commutators in interpolation theory}, Analysis at Urbana, Vol. II (Urbana, IL, 1986--1987), 170--220, London Math. Soc. Lecture Note Ser., 138. Cambridge Univ. Press, Cambridge, 1989. 


\bibitem{CMR} M. Cwikel, M. Milman and R. Rochberg, {\it An introduction to Nigel Kalton's work on differentials of complex interpolation processes for K\"othe spaces}, arXiv:1404.2893.


\bibitem{DJT} J. Diestel, H. Jarchow and A. Tongue, \emph{Absolutely summing operators}. Cambridge Studies in Advanced Mathematics, 43. Cambridge University Press, Cambridge, 1995.

\bibitem{CoS} F. Cobos and T. Schonbek, {\it On a theorem by Lions and Peetre about interpolation between a Banach space and its dual}, Houston J. Math. 24 (1998), no. 2, 325--344.

\bibitem{CJM} M. Cwikel, B. Jawerth and M. Milman, {\it The domain spaces of quasilogarithmic operators}, Trans. Amer. Math. Soc. 317 (1990), no. 2, 599--609. 

\bibitem{JRW} B. Jawerth, R. Rochberg and G. Weiss, {\it Commutator and other second order estimates in real interpolation theory}, Ark. Mat. 24 (1985), 191--219.

\bibitem{JoSa} W.B. Johnson and A. Szankowski, {\it Hereditary approximation property}, Ann. of Math. (2) 176 (2012), no. 3, 1987--2001.

\bibitem{KaPe} N.J. Kalton and N.T. Peck, {\it Twisted sums of sequence spaces and the three-space problem}, Trans. Amer. Math. Soc. 255 (1979), 1--30.

\bibitem{KM} N.J. Kalton and S. Montgomery-Smith, {\it Interpolation of Banach spaces}, Handbook of Geometry of Banach Spaces, Vol. 2, (W.B. Johnson and J. Lindenstrauss, editors), Elsevier, Amsterdam, 2003, 1131--1175.

 \bibitem{ka} N.J. Kalton, {\it Differentials of complex interpolation processes for K\"othe function spaces}, Trans. Amer. Math. Soc. 333 (1992), no. 2, 479--529.

%
%



%
%
%
%
%
%
%
%
%
\bibitem{K4} N.J. Kalton, {\it The three space problem for locally bounded F-spaces}, Compo. Math. 37 (1978), 243--276.

%
\bibitem{LT} J. Lindenstrauss and L. Tzafriri, \emph{Classical Banach
spaces I, sequence spaces}, Ergebnisse der Math. und ihrer
Grenzgebiete 92, Springer-Verlag 1977.

\bibitem{LiPe} J.L. Lions and J. Peetre, {\it Sur une classe d'espaces d'interpolation},  Inst. Hautes \'Etudes Sci. Publ. Math. (1964), No. 19, 5--68.

\bibitem{JS3} D. Morales and J. Su\'arez de la Fuente, {\it Some more twisted Hilbert spaces}, Ann. Fenn. Math. 46 (2021), no. 2, 819--837.

\bibitem{Pi} G. Pisier, {\it Weak Hilbert spaces}, Proc. London Math. Soc. (3) 56 (1988), no. 3, 547--579.

\bibitem{Pib} G. Pisier, \emph{The volume of convex bodies and Banach space geometry}. Cambridge Tracts in Mathematics, 94. Cambridge University Press, Cambridge, 1989.


\bibitem{RW} R. Rochberg and G. Weiss, {\it Derivatives of Analytic Families of Banach Spaces}, Annals of Mathematics 118, no. 2 (1983), 315--347.

\bibitem{Se} S. Semmes, {\it Interpolation of Banach spaces, differential geometry and differential equations}, Rev. Mat. Iberoamericana (1988), no. 1, 155--176.

\bibitem{JS} J. Su\'arez de la Fuente, {\it A weak Hilbert space that is a twisted Hilbert space}, J. Inst. Math. Jussieu 19 (2020), no. 3, 855--867.

\bibitem{JS2} J. Su\'arez de la Fuente, {\it A space with no unconditional basis that satisfies the Johnson-Lindenstrauss lemma},  Results Math. 74 (2019), no. 3, Art. 126, 14 pp.



\end{thebibliography}
\end{document}